\documentclass[twoside,leqno]{article}
\usepackage{latexsym,amsmath,amssymb}

\catcode`\@=11
\def\keywords#1{\def\@keywords{#1}}
\let\@keywords=\@empty
\def\subjclass#1{\def\@subjclass{#1}}
\let\@subjclass=\@empty
\newcommand{\keywordsname}{Key words and phrases}
\newcommand{\subjclassname}{\textup{1991} Mathematics Subject
     Classification}
\def\@addpunct#1{\ifnum\spacefactor>\@m \else#1\fi}
\def\@setsubjclass{%
  {\itshape\subjclassname.}\enspace\@subjclass\@addpunct.}
\def\@setkeywords{%
  {\itshape \keywordsname.}\enspace \@keywords\@addpunct.}
\renewcommand\thesection       {\arabic{section}}

\def\initrhead{}%
\def\initrhead{{\it JAMS}}%
\newcounter{footpage}%
\def\sikibetuh{}%
\mark{{}{}}%
\def\ps@myheadings{%
    \let\@oddfoot\@empty\let\@evenfoot\@empty
    \def\@evenhead{\small\sikibetuh\thepage\hfill\leftmark\hfill}%
    \def\@oddhead{\small\hfill\rightmark\hfill\sikibetuh\thepage}%
    \let\@mkboth\@gobbletwo
    \let\sectionmark\@gobble
    \let\subsectionmark\@gobble
    }
\def\ps@jamsinit{%
    \def\@evenhead{\small\sikibetuh\thepage\hfill{\footnotesize{\initrhead}}}%
    \def\@oddhead{\small{\footnotesize{\initrhead}}\hfill\sikibetuh\thepage}%
    \let\@mkboth\@gobbletwo
    \let\sectionmark\@gobble
    \let\subsectionmark\@gobble
    }
\def\section{\@startsection{section}{1}%
  \z@{.7\baselineskip\@plus\baselineskip}{-.5em}%
  {\normalfont\bfseries}}
\def\subsection{\@startsection{subsection}{2}%
  \z@{.5\baselineskip\@plus.7\baselineskip}{-.5em}%
  {\normalfont\bfseries}}
\def\subsubsection{\@startsection{subsubsection}{3}%
  \z@{.5\baselineskip\@plus.7\baselineskip}{-.5em}%
  {\normalfont\itshape}}
\def\paragraph{\@startsection{paragraph}{4}%
  \z@\z@{-\fontdimen2\font}%
  \normalfont}
\def\subparagraph{\@startsection{subparagraph}{5}%
  \z@\z@{-\fontdimen2\font}%
  \normalfont}
\def\appendix{\par\c@section\z@ \c@subsection\z@
   \let\sectionname\appendixname
   \def\thesection{\@Alph\c@section}}
\def\appendixname{Appendix}
\renewenvironment{abstract}{%
\noindent
  \small
\noindent
  \quote
  \parindent=1.5em
  {\scshape \abstractname .\hspace{\z@}}%
  }
{  \endquote
}
\renewcommand{\maketitle}{\par
  \begingroup
    \renewcommand{\thefootnote}{\fnsymbol{footnote}}%
    \def\@makefnmark{\hbox to\z@{$\m@th^{\@thefnmark}$\hss}}%
    \long\def\@makefntext##1{\parindent 1em\noindent
            \hbox to1.8em{\hss$\m@th^{\@thefnmark}$}##1}%
      \newpage
      \global\@topnum\z@   %
      \@maketitle
    \@thanks
  \endgroup
  \thispagestyle{jamsinit}
  \setcounter{footnote}{0}%
  \let\thanks\relax
  \let\maketitle\relax\let\@maketitle\relax
  \gdef\@thanks{}
  \gdef\@author{}\gdef\@title{}}
\def\@maketitle{%
  \newpage
  \null
  \vskip 2em%
  \begin{center}%
    {\normalfont\bfseries \@title \par}%
    \vskip 1.5em%
    {\normalfont\scshape
      \lineskip .5em%
      \begin{tabular}[t]{c}%
        \@author
      \end{tabular}\par
     \def\thefootnote{}
     \ifx\@empty\@subjclass\else \footnotetext{\@setsubjclass}\fi
     \def\thefootnote{}
     \ifx\@empty\@keywords\else \footnotetext{\@setkeywords}\fi
}%
    \vskip 1em%
    {\normalfont \@date}%
  \end{center}%
  \par
  \vskip 0.8em
}
\renewenvironment{thebibliography}[1]
      {
  \@startsection{section}
        \@m\z@{9\p@\@plus12\p@}{6\p@}%
        {\centering\scshape}\refname
  \small\labelsep .5em\relax
      \list{\@biblabel{\@arabic\c@enumiv}}%
           {\settowidth\labelwidth{\@biblabel{#1}}%
            \leftmargin\labelwidth
            \advance\leftmargin\labelsep
            \@openbib@code
            \usecounter{enumiv}%
            \let\p@enumiv\@empty
            \renewcommand\theenumiv{\@arabic\c@enumiv}}%
      \sloppy
      \clubpenalty4000
      \@clubpenalty \clubpenalty
      \widowpenalty4000%
      \sfcode`\.\@m}
     {\def\@noitemerr
       {\@latex@warning{Empty `thebibliography' environment}}%
      \endlist}
\newtheorem{theorem}{Theorem}
\newtheorem{corollary}[theorem]{Corollary}
\newtheorem{definition}[theorem]{Definition}
\newtheorem{example}[theorem]{Example}
\newtheorem{lemma}[theorem]{Lemma}
\newtheorem{proposition}[theorem]{Proposition}
\newtheorem{remark}[theorem]{Remark}
\begin{document}
\title{Covering theorems and Lebesgue integration}
\author{Peter A. Loeb\thanks{Supported in part by NSF Grant DMS96-22454.}\\Department of Mathematics, University of 
Illinois\\1409 West Green Street, 
Urbana, Illinois 61801, U.S.A.\\e-mail: loeb@math.uiuc.edu
\and Erik Talvila\thanks{Supported by an NSERC Postdoctoral
Fellowship.}\\Department of Mathematical Sciences, University of 
Alberta\\Edmonton, 
Alberta T6G 2E2, Canada\\e-mail: etalvila@math.ualberta.ca}
\date{June 13, 2000}
\subjclass{Primary 28A25, 28A75, 52C17}
\keywords{Besicovitch Covering Theorem, Morse Covering Theorem, Open Covering
Theorem, Lebesgue integral, Riemann sum}
\maketitle
\begin{abstract}
This paper shows how the Lebesgue integral can be obtained as a Riemann sum
and provides an extension of the Morse Covering Theorem to open sets. Let $X$
be a finite dimensional normed space; let $\mu$ be a Radon measure on $X$ and
let $\Omega\subseteq X$ be a $\mu$-measurable set. For $\lambda\geq1$, a $\mu
$-measurable set $S_{\lambda}(a)\subseteq X$ is a $\lambda$-Morse set with tag
$a\in S_{\lambda}(a)$ if there is $r>0$ such that $B(a,r)\subseteq S_{\lambda
}(a)\subseteq B(a,\lambda r)$ and $S_{\lambda}(a)$ is starlike with respect to
all points in the closed ball $B(a,r)$. Given a gauge $\delta\!:\!\Omega
\rightarrow(0,1]$ we say $S_{\lambda}(a)$ is $\delta$-fine if $B(a,\lambda
r)\subseteq B(a,\delta(a))$. If $f\geq0$ is a $\mu$-measurable function on
$\Omega$ then $\int_{\Omega}f\,d\mu=F\in\mathbb{R}$ if and only if for some
$\lambda\geq1$ and all $\varepsilon>0$ there is a gauge function $\delta$ so
that $|\sum_{n}f(x_{n})\,\mu(S(x_{n}))-F|<\varepsilon$ for all sequences of
disjoint $\lambda$-Morse sets that are $\delta$-fine and cover all but a $\mu
$-null subset of $\Omega$. This procedure can be applied separately to the
positive and negative parts of a real-valued function on $\Omega$. The
covering condition $\mu(\Omega\setminus\cup_{n}S(x_{n}))=0$ can be satisfied
due to the Morse Covering Theorem. The improved version given here says that
for a fixed $\lambda\geq1$, if $A$ is the set of centers of a family of
$\lambda$-Morse sets then $A$ can be covered with the interiors of sets from
at most $\kappa$ pairwise disjoint subfamilies of the original family; an
estimate for $\kappa$ is given in terms of $\lambda$, $X$ and its norm.
\end{abstract}

\section{Introduction}

An attractive feature of the Riemann and Henstock integrals is that they can
be defined in terms of Riemann sums. Suppose we wish to integrate a
real-valued function $f$ over a set $\Omega$ with respect to a measure $\mu$.
If we have disjoint measurable sets $\Omega_{1},\ldots,\Omega_{N}$ with union
$\Omega$ (i.e., a partition of $\Omega$), then we may try to define an
integral as the limit of sums $\sum_{i=1}^{N}\!f(z_{i})\,\mu(\Omega_{i})$ for
appropriate points $z_{i}\in\Omega$. One would hope that taking the sets
$\Omega_{i}$ small enough and $N$ large enough would make these sums close to
the same value, which we then define to be the integral $\int_{\Omega}f\,d\mu
$. When $\Omega\subseteq\mathbb{R}^{d}$, this is done in the Riemann case for
Lebesgue measure and a bounded function $f$ and bounded set $\Omega$ by
choosing for the partition sets $\Omega_{i}$ uniformly small cubes and then
choosing arbitrary points $z_{i}\in\Omega_{i}$. With the Henstock integral,
$f$ and $\Omega$ need no longer be bounded. For this case, the sets
$\Omega_{i}$ are intervals satisfying a gauge condition. This means to begin
with that we have a function $\delta\!:\!\Omega\rightarrow(0,R)$ for some
positive $R$; the mapping is called a \textbf{gauge function}, and we say the
pair $(z_{i},\Omega_{i})$ is $\delta$\textbf{-fine} if $z_{i}\in\Omega_{i}$
and $\Omega_{i}$ is contained in the closed ball with center $z_{i}$ and
radius $\delta(z_{i})$. (When $\Omega$ is unbounded, the partition need
only cover $\Omega\cap B(0,R)$ where $B(0,R)$ is a ball with center at
the origin and large enough radius $R$ determined by 
$\delta$.) We obtain the McShane
integral by dropping the restriction that $z_{i}\in\Omega_{i}$. See
\cite{gordon} for a discussion of these integrals.

All of these integration schemes revolve around finding a partition of
$\Omega$, which of course requires rather specialized sets $\Omega_{i}$.
Breaking this pattern, the Vitali Covering Theorem was used in \cite{malee} to
define the Lebesgue integral with respect to Lebesgue measure on a finite
interval of the real line. There the idea, given any $\eta>0$, is to use a
finite number of intervals so that $\lambda([a,b]\setminus\cup_{i=1}^{N}
I_{i})<\eta$. Here, we too will apply covering theory, but now to a measurable
set $\Omega$ in a finite dimensional normed space $X$. We will obtain the
Lebesgue integral with respect to a Radon measure as a series $\sum
f(z_{i})\mu(\Omega_{i})$, where the sets $\Omega_{i}\subseteq X$ are disjoint
and cover all but a null set of $\Omega$. The sets $\Omega_{i}$ will be made
small with respect to a gauge function. They can be balls or starlike sets
(described in Section~\ref{Covering} below). It is the Besicovitch Covering
Theorem for balls and the Morse Covering Theorem for starlike sets that
enables us to fulfill the condition $\mu(\Omega\setminus\cup_{i}\Omega_{i}
)=0$. For this it is essential that we are working in a finite dimensional
normed space, not just a metric space. These covering results are discussed
below, and simplified proofs of strengthened versions are provided. We also
note that by omitting a small part of the overall sum $\sum f(z_{i})\mu
(\Omega_{i})$, we are able to restrict the points $z_{i}$ to the set of points
of approximate continuity of $f$, defined in terms of Morse covers in
Section~\ref{measure}.

In the theory of Henstock and McShane integration, the appearance of the gauge
function is rather mysterious: For all $\varepsilon>0$ there is a gauge
function $\delta\!:\!\Omega\rightarrow(0,\infty)$ such that for all $\delta
$-fine partitions $\{(z_{i},\Omega_{i})\}_{i=1}^{N}$ of $\Omega$ we have
$|\sum_{i=1}^{N}\!f(z_{i})\,\mu(\Omega_{i})-\int_{\Omega}f\,d\mu|<\varepsilon
$. We show in proving Theorem \ref{integration} how the properties of Lebesgue
points can be used to determine the gauge $\delta$. An even simpler result
extending the Riemann integral is obtained in Section \ref{Riemann} for the
case that $f$ is continuous at $\mu$-almost all points of $\Omega$.

\section{Covering Theorems}

\label{Covering}

Let $\left(  X,\left\|  \cdot\right\|  \right)  $ be a normed vector space of
dimension $d<\infty$ over the real numbers $\mathbb{R}$. Then $X$ is a
separable, locally compact Hausdorff space with open sets determined by the
open balls. The open ball with center $\mathbf{a}\in X$ and radius $r>0$ is
denoted by $U(\mathbf{a},r):=\{\mathbf{x}\in X:\Vert\mathbf{x}-\mathbf{a}
\Vert<r\}$; the closed ball with center $\mathbf{a}\in X$ and radius $r>0$ is
$B(\mathbf{a},r):=\{\mathbf{x}\in X:\Vert\mathbf{x}-\mathbf{a}\Vert\leq r\}$.
Since the dimension of $X$ is finite, its closed balls are compact.

The integration results to follow will use coverings by balls and sets more
general than balls in $\left(  X,\left\|  \cdot\right\|  \right)  $, and for
this we will need the Besicovitch and Morse covering theorems. Strengthened
versions of these theorems are as easy to state and prove as the original
results. This was done at a real-analysis meeting in Rolla, Missouri in 1995,
and the work was included in the report of that meeting as the note by the
first author in \cite{loeb}. Since all of that report is now essentially
unavailable, we will sketch these improved results and proofs here for the
reader's convenience.

For general finite dimensional normed vector spaces, the Besicovitch Covering
Theorem uses covers by closed metric balls; it gives a constant that is
independent of measure. Besicovitch's result is much stronger than the
familiar Vitali Covering Theorem. It was originally established for disks in
the plane in 1945-46 \cite{besicovitch}, and was extended by A. P. Morse
\cite{morse} in 1947 to more general shapes in finite dimensional normed
spaces. The constructions used in both the Besicovitch and Morse results are
modified here so that better theorems are obtained. In the modified theorems,
the original cover of a set $A$ by closed sets can still be reduced to a
subcover $\mathcal{F}$ such that $\mathcal{F}$ can be partitioned into $n$
subfamilies of pairwise disjoint sets and $n$ is bounded above by a global
constant depending only on the space. The construction of $\mathcal{F}$
\thinspace is arranged, however, so that $A$ is actually contained in the
union of the interiors of the sets in $\mathcal{F}$. To obtain this result, we
have modified the following definition taken from \cite{bliedtnerloeb1}.

\begin{definition}
Fix $\tau>1$. Let $\{S_{i}:1\leq i\leq n\}$ be an ordered collection of
subsets of $X$ with each $S_{i}$ having finite diameter $\Delta(S_{i})$ and
containing a point $a_{i}$ in its interior\textbf{, }$\operatorname*{int}
(S_{i})$. We say that the ordered collection of sets $S_{i}$ is in $\tau
$-satellite configuration with respect to the ordered set of points $a_{i}$ if
\textbf{i)}~For all $i\leq n$,\ $S_{i}\cap S_{n}\neq\varnothing$ and
\textbf{ii)} For all pairs $i<j\leq n$, \ $a_{j}\notin\operatorname*{int}
(S{}_{i})$ \thinspace and$\ \,\,\Delta(S_{j})<\tau\cdot\Delta(S_{i})$.

\begin{theorem}
\label{Cover}Let $A$ be an arbitrary subset of $X$. With each point $a\in A$,
associate a set $S(a)$ containing $a$ in its interior so that the diameters
have a finite upper bound. Assume that for some $\tau>1$, there is an upper
bound $\kappa\in\mathbb{N}\mathbf{\ }$ to the cardinality of any ordered set
$\{a_{i}:1\leq i\leq n\}\subseteq A$ with respect to which the ordered set
$\{S(a_{i}):1\leq i\leq n\}$ is in $\tau$-satellite configuration. Then for
some $m\leq\kappa$, there are pairwise disjoint subsets $A_{1},\ldots,A_{m}$
of $A$ such that $A\subseteq\cup_{j=1}^{m}\cup_{a\in 
A_{j}}\operatorname*{int}
(S\,(a))$ and for each $j$, $1\leq j\leq m$, the elements of the collection
$\{S(a):a\in A_{j}\}$ are pairwise disjoint.
\end{theorem}
\end{definition}

\textbf{Proof:} Let $T$ be a choice function on the nonempty subsets $B$ of
$A$ such that $T(B)$ is a point $b\in B$ with $\tau\cdot\Delta(S(b))>\sup
_{a\in B}\Delta(S(a))$. Form a one-to-one correspondence between an initial
segment of the ordinal numbers and a subcollection of $A$ as follows. Set
$B_{1}=A$ and $a_{1}=T(B_{1})$. Having chosen $a_{\alpha}$ for $\alpha<\beta$,
let $B_{\beta}=A\setminus\cup_{\alpha<\beta}\operatorname*{int}(S{}(a_{\alpha
}))$. If $B_{\beta}\neq\varnothing$, set $a_{\beta}=T(B_{\beta})$. There
exists a first ordinal $\gamma$ for which $B_{\gamma}=\varnothing$; that is,
$A\subseteq\cup_{\alpha<\gamma}\operatorname*{int}(S{}(a_{\alpha}))$. Note
that for $\alpha<\beta<\gamma$, we have $a_{\beta}\notin\operatorname*{int}
(S{}(a_{\alpha}))\;$and$\ \ \Delta(S(a_{\beta}))<\tau\cdot\Delta(S(a_{\alpha
}))$. Let $A_{c}=\{a_{\alpha}:\alpha<\gamma\}$, and let $\prec$ denote the
well-ordering on $A_{c}$ inherited from the ordinals.

Given any nonempty subset $B$ of $A_{c}$, form a one-to-one correspondence
between an initial segment of the ordinal numbers and a subset $V(B)$ of $B$
as follows. Set $B_{1}=B$, and let $a(1)$ be the first element (with respect
to $\prec$) of $B_{1}$. Having chosen $a(\alpha)$ for $\alpha<\beta$, let
\[
B_{\beta}=\{b\in B:\;\forall\alpha<\beta,\;S(b)\cap S(a(\alpha))=\varnothing
\}.
\]
If $B_{\beta}\neq\varnothing$, let $a(\beta)$ equal the first element (with
respect to $\prec$) of $B_{\beta}$. There exists a first ordinal $\gamma$ for
which $B_{\gamma}=\varnothing$. Let $V(B)=\{a(\alpha):\alpha<\gamma\}$.

Now for $i\geq1$, form sets $A_{i}\subseteq A_{c}$ as follows. Set
$A_{1}=V(A_{c})$. Having chosen $A_{i}$ for $1\leq i\leq n$, let $B_{n}
=A_{c}\setminus\cup_{i=1}^{n}A_{i}$. Stop if $B_{n}=\varnothing$. Otherwise,
set $A_{n+1}=V(B_{n})$. Note that for each $b\in B_{n}$ and each $i$ between
$1$ and $n$, there is a first (with respect to $\prec$) $a_{i}\in A_{i}$ with
$S(a_{i})\cap S(b)\neq\varnothing$; clearly, $a_{i}\prec b$ in $A_{c}$. It now
follows that the set $\{S(a_{1}),\ldots,S(a_{n}),S(b)\}$ is in $\tau
$-satellite configuration with respect to the set $\{a_{1},\ldots,a_{n},b\}$
when each set is given the ordering inherited from $A_{c}$. Therefore,
$B_{n}=\varnothing$ for some $n\leq\kappa$. $\Box$

\begin{corollary}
\label{finiteborelmeasure} For any finite Borel measure $\mu$ on $X$, there is
a $j$ with $1\leq j\leq m$ and a finite subset $\;A_{\mu}\subseteq A_{j}\;$
such that
\[
\mu^{\ast}(A)\leq2\kappa\cdot\sum_{a\in A_{\mu}}\mu(\operatorname*{int}
(S\,(a))).
\]
\end{corollary}

\bigskip\textbf{Proof:} Take the first $j\leq m$ that maximizes the sum
$\Sigma_{a\in A_{j}}\,\mu(\operatorname*{int}(S\,(a)))$. We can then choose a
finite subset $A_{\mu}\subseteq A_{j}$ so that $\frac{1}{2}\cdot\sum_{a\in
A_{j}}\mu(\operatorname*{int}(S\,(a)))\leq\sum_{a\in A_{\mu}}\mu
(\operatorname*{int}(S\,(a)))$.\ $\Box$

\bigskip

What is the upper bound $\kappa$ for our vector space $X$? For balls and
values of $\tau$ close to $1$, there is an upper bound $K$ for $\kappa$
established by Zolt\'{a}n F\"{u}redi and the first author in 
\cite{furediloeb}.
It is the maximum number of points that can be packed into the closed ball
$B(\mathbf{0},2)$ when one of the points is at $\mathbf{0}$ and the distance
between distinct points is at least $1$. That value is no more than $5^{d}$,
where $d$ is the dimension of $X$. Applying Theorem~\ref{Cover} with
$\kappa=K$ \ to a cover by balls yields an open version of Besicovitch's
theorem for $X$. The constant $K$ is the best constant for the Besicovitch
Theorem in terms of all known proofs. With obvious modifications, the
construction in \cite{furediloeb} is already appropriate for the improved
result that yields a cover by open balls.

To use Theorem~\ref{Cover} to establish an open version of Morse's Covering
Theorem for $(X,\left\|  \cdot\right\|  )$, we need some geometric results.
The proofs are modifications of arguments in \cite{morse} and
\cite{bliedtnerloeb1}. The bound obtained is not as simple as the one for
balls, but the shapes to which it applies are more general than balls or even
convex sets. For these geometric arguments, we use boldface to denote points.

For each $\gamma\geq1$, we let $N(\gamma)$ be an upper bound for the number of
points that can be packed into the closed ball $B(\mathbf{0},1)$ when the
distance between distinct points is at least $1/\gamma$ and one point is at
$\mathbf{0}$. We write $N_{S}(\gamma)$ for the similar constant when all
points are on the surface of $B(\mathbf{0},1)$. Given nonzero points
$\mathbf{b}$ and $\mathbf{c}$ in $X$, we set $V(\mathbf{b},\mathbf{c}
):=\left\|  \frac{\mathbf{b}}{\left\|  \mathbf{b}\right\|  }-\frac{\mathbf{c}
}{\left\|  \mathbf{c}\right\|  }\right\|  $.

\begin{proposition}
\label{Morse1}Fix $\tau$ with $1<\tau\leq2$. Also fix an ordered set
$\{S_{i}:1\leq i\leq n\}$ of bounded subsets of $X$ each containing a ball
$B(\mathbf{a}_{i},r_{i})$. Assume that $\{S_{i}:1\leq i\leq n\}$ is in $\tau
$-satellite configuration with respect to the ordered set of centers
$\{\mathbf{a}_{i}:1\leq i\leq n\}$. Translate so that $\mathbf{a}
_{n}=\mathbf{0}$. Fix $\lambda\geq\max_{1\leq i\leq n}\Delta(S_{i})/(2r_{i})$.
Suppose the resulting configuration has the following property in terms of two
constants $C_{0}\geq1$ and $C_{1}\geq1:$ If $\mathbf{a}_{i}$ and
$\mathbf{a}_{j}$ are centers with the properties that $C_{0}r_n<\Vert
\mathbf{a}_{i}\Vert\leq\Vert\mathbf{a}_{j}\Vert$ \ \ and $\ V(\mathbf{a}
_{i},\mathbf{a}_{j})\leq1/C_{1}$, then $\mathbf{a}_{i}$ must be in the
interior of $S_{j}$. It then follows that
\[
n\leq N(2\lambda C_{0})+N(8\lambda^{2})\ N_{S}(C_{1}).
\]
\end{proposition}

\textbf{Proof:} Set \thinspace$r:=r_{n}$ and \thinspace$S:=S_{n}$. For $1\leq
i<j\leq n$, we have
\[
\Vert\mathbf{a}_{i}-\mathbf{a}_{j}\Vert\geq r_{i}\geq\Delta(S_{i}
)/(2\lambda)\geq\Delta(S)/(4\lambda)\geq r/(2\lambda).
\]
Scaling by $1/(C_{0}r)$, one sees that there can be at most $N(2\lambda
C_{0})$ indices $i$ for which $\Vert\mathbf{a}_{i}\Vert\leq C_{0}r$. We only
have to show, therefore, that there are at most $N(8\lambda^{2})\ N_{S}
(C_{1})$ indices in the set $J:=\{j<n:C_{0}r<\Vert\mathbf{a}_{j}\Vert\}$.
Suppose $i\neq j$ are members of $J$ with $\mathbf{a}_{i}\in
\operatorname*{int}(S{}_{j})$. Then $i<j$ and
\[
\mathbf{a}_{j}\in B(\mathbf{a}_{i},\Delta(S_{j}))\subseteq B(\mathbf{a}
_{i},2\Delta(S_{i}))\subseteq B(\mathbf{a}_{i},4\lambda r_{i}).
\]
Moreover, $\Vert\mathbf{a}_{j}-\mathbf{a}_{i}\Vert\geq r_{i}\geq
r_{i}/(2\lambda)$. If also $j<k$ in $J$, and $\mathbf{a}_{i}\in
\operatorname*{int}(S{}_{k})$, then\ $\mathbf{a}_{k}\in B(\mathbf{a}
_{i},4\lambda r_{i})$ \ and
\[
\Vert\mathbf{a}_{k}-\mathbf{a}_{j}\Vert\geq r_{j}\geq\Delta(S_{j}
)/(2\lambda)\geq\Vert\mathbf{a}_{j}-\mathbf{a}_{i}\Vert/(2\lambda)\geq
r_{i}/(2\lambda)\text{.}
\]
Scaling by $1/(4\lambda r_{i})$, it follows that for each $i\in J$, the
cardinality $\operatorname*{Card}\{j\in 
J:\mathbf{a}_{i}\in\operatorname*{int}
(S{}_{j})\}\leq N(8\lambda^{2})$. Now construct $J^{\prime}\subseteq J$ by
induction as follows. Set $J_{1}=J$. At the $k^{\text{th}}$ step for $k\geq1$,
if $J_{k}$ is empty, stop. Otherwise, choose the first $i_{k}\in J_{k}$ so
that for all $j\in J_{k}$, $\Vert\mathbf{a}_{i_{k}}\Vert\leq\Vert
\mathbf{a}_{j}\Vert$. Put $i_{k}$ in $J^{\prime}$. Form the set $J_{k+1}$ by
discarding from $J_{k}$ the index $i_{k}$ and all other indices $j$ such that
$\mathbf{a}_{i_{k}}\in\operatorname*{int}(S{}_{j})$. Now, if $i\neq j$ in
$J^{\prime}$, $V(\mathbf{a}_{i},\mathbf{a}_{j})>1/C_{1}$. Therefore,
$\operatorname*{Card}(J^{\prime})\leq N_{S}(C_{1})$, and so
$\operatorname*{Card}(J)\leq N(8\lambda^{2})\ N_{S}(C_{1})$. $\Box$

\bigskip

Given $\lambda\geq1$ and $\mathbf{a}\in X$, we let $\,\mathcal{S}_{\lambda
}(\mathbf{a})\,$ denote the collection of all sets $S\subseteq X$ \thinspace
for which there exists an $r>0$ such that $B(\mathbf{a},r)\subseteq S\subseteq
B(\mathbf{a},\lambda r)$ and $S$ is \textbf{starlike} with respect to every
$\mathbf{y}\in B(\mathbf{a},r)$. This means that for each $\mathbf{y}\in
B(\mathbf{a},r)$ and each $\mathbf{x}\in S$, the line segment $\alpha
\mathbf{y}+(1-\alpha)\mathbf{x}$, $0\leq\alpha\leq1$, is contained in $S$.
This is the general shape considered by Morse in \cite{morse}. To improve his
result, as well as for work in a later section, we will need the following
fact about such a set $S$; the result, along with the next theorem, will
finish our proof of the ``open'' Morse's Covering Theorem.

\begin{proposition}
\label{cone}If $\left\|  \mathbf{y}-\mathbf{a}\right\|  <r$, i.e., if
$\mathbf{y}$ is in the interior of $B(\mathbf{a},r)$, and $\mathbf{x}$ is in
the closure, $\operatorname*{cl}(S)$, of $S$, then every point of the form
$\alpha\mathbf{y}+(1-\alpha)\mathbf{x}$, $0<\alpha\leq1$, is in the interior
of $S$.
\end{proposition}

\textbf{Proof:} Fix $\rho>0$ so that $B(\mathbf{y},\rho)\subset B(\mathbf{a}
,r)$, and fix $\alpha$ with $0<\alpha\leq1$. Assume first that $\mathbf{x}\in
S$, and translate so that $\mathbf{x}=\mathbf{0}$. Then the ball
$B(\alpha\mathbf{y},\alpha\rho)\subseteq S$ since
\begin{align*}
\left\|  \alpha\mathbf{y}-\mathbf{z}\right\|   &  \leq\alpha\rho
\Rightarrow\left\|  \mathbf{y}-\tfrac{1}{\alpha}\mathbf{z}\right\|  \leq
\rho\Rightarrow\tfrac{1}{\alpha}\mathbf{z}\in B(\mathbf{a},r)\\
&  \Rightarrow\mathbf{z}=\alpha\left(  \tfrac{1}{\alpha}\mathbf{z}\right)
+(1-\alpha)\mathbf{0}\in S\text{.}
\end{align*}
Now for the case that $\mathbf{x}\in\operatorname*{cl}(S)$, choose a point
$\mathbf{w}\in S$ so that $\frac{1-\alpha}{\alpha}\left\|  \mathbf{x}
-\mathbf{w}\right\|  <\rho$. The result follows from the previous case since
\[
\alpha\mathbf{y}+(1-\alpha)\mathbf{x}=\alpha\left(  \mathbf{y}+\tfrac
{1-\alpha}{\alpha}\left(  \mathbf{x}-\mathbf{w}\right)  \right)
+(1-\alpha)\mathbf{w}\text{. \ }\Box
\]

\begin{theorem}
Fix $\lambda\geq1$ and fix $\tau$ with $1<\tau\leq2$. If $\{S_{i}:1\leq i\leq
n\}$ is an ordered collection of subsets of $X$ in $\tau$-satellite
configuration with respect to an ordered set $\{\mathbf{a}_{i}:1\leq i\leq
n\}\subset X$, and if for $1\leq i\leq n$, $S_{i}\in\mathcal{S}_{\lambda
}(\mathbf{a}_{i})$, then
\[
n\leq N(64\lambda^{3})+N(8\lambda^{2})N_{S}(16\lambda).
\]
\end{theorem}

\textbf{Proof:} For $1\leq i\leq n$, fix $r_{i}>0$ so that $B(\mathbf{a}
_{i},r_{i})\subseteq S_{i}\subseteq B(\mathbf{a}_{i},\lambda r_{i})$ and
$S_{i}$ is starlike with respect to every $\mathbf{y}\in B(\mathbf{a}
_{i},r_{i})$. Translate so that $\mathbf{a}_{n}=\mathbf{0}$; set $r=r_{n}$ and
$S=S_{n}$. Suppose $i$ and $j$ are indices such that $32\lambda^{2}
r<\Vert\mathbf{a}_{i}\Vert\leq\Vert\mathbf{a}_{j}\Vert$ and $V(\mathbf{a}
_{i},\mathbf{a}_{j})\leq1/(16\lambda)$. By Proposition \ref{Morse1}, we only
have to show that $\mathbf{a}_{i}$ must be in the interior of $S{}_{j}$. To
simplify notation, let $\mathbf{b}=\mathbf{a}_{i}$ and $\mathbf{c}
=\mathbf{a}_{j}$. Fix $\mathbf{x}\in S\cap S_{j}$. Since $\left\|
\mathbf{x}\right\|  \leq\lambda r<32\lambda^{2}r<\Vert\mathbf{b}\Vert$,
$\;\mathbf{x}\neq\mathbf{b}$. Let $s=\Vert\mathbf{c}\Vert/\Vert\mathbf{b}
\Vert$ and $t=1/s$. Set $\mathbf{y}=(1-s)\mathbf{x}+s\mathbf{b}$. Then
$\mathbf{b}=(1-t)\mathbf{x}+t\mathbf{y}$. To show that $\mathbf{b}
\in\operatorname*{int}(S{}_{j})$, we only have to show that$\;\Vert
\mathbf{y}-\mathbf{c}\Vert<r_{j}$.\ Now $16\lambda\Delta(S)\leq32\lambda
^{2}r<\Vert\mathbf{b}\Vert$, whence $\Vert\mathbf{x}\Vert\leq\Delta(S)\leq
\min\left(  \Vert\mathbf{b}\Vert/(16\lambda),\,2\Delta(S_{j})\right)  $.
Therefore, since $\left|  1-s\right|  =s-1<s$, $\;$
\begin{align*}
\Vert\mathbf{y}-\mathbf{c}\Vert &  =\;\left\|  (1-s)\mathbf{x}+\left\|
\mathbf{c}\right\|  \left(  \tfrac{\mathbf{b}}{\Vert\mathbf{b}\Vert}
-\tfrac{\mathbf{c}}{\Vert\mathbf{c}\Vert}\right)  \right\| \\
&  <s\Vert\mathbf{x}\Vert+\Vert\mathbf{c}\Vert/(16\lambda)\\
&  \leq s\Vert\mathbf{b}\Vert/(16\lambda)+\Vert\mathbf{c}\Vert/(16\lambda
)=\Vert\mathbf{c}\Vert/(8\lambda)\\
&  \leq\left(  \Vert\mathbf{c}-\mathbf{x}\Vert+\Vert\mathbf{x}\Vert\right)
/(8\lambda)\\
&  <\Delta(S_{j})/(2\lambda)\leq r_{j}.\;\;\Box
\end{align*}

\section{Measures}

\label{measure}

Recall that we are working with a normed vector space $(X,\Vert\cdot\Vert)$ of
dimension $d<\infty$ over the real numbers $\mathbb{R}$. Let $\mu$ be a
measure on a $\sigma$-algebra $\mathcal{M}$ of subsets of $X$. We say that
$\mu$ is a \textbf{Radon measure} on $X$ if:

\noindent\textbf{(i)}\ All Borel sets are measurable, i.e., $\mathcal{M}$
contains the Borel sets.

\noindent\textbf{(ii)}\ Compact sets have finite measure.

\noindent\textbf{(iii)} $\mu$ is \textbf{inner and outer regular}, i.e., for
all $E\in\mathcal{M}$
\begin{align*}
\mu(E)  &  =\sup\{\mu(K):K\subseteq E\text{ and }K\text{ is compact}\},\\
\mu(E)  &  =\inf\{\mu(G):G\supseteq E\text{ and }G\text{ is open}\}.
\end{align*}
We will call a set or function $\mu$\textbf{-measurable}, or when $\mu$ is
understood just \textbf{measurable}, if it is measurable with respect to the
$\mu$-completion of $\mathcal{M}$.

Since every open set in $X$ is $\sigma$-compact, inner and outer regularity
follow from assuming merely that $\mu$ is a Borel measure on $X$ and closed
balls have finite measure; see Theorem~2.18 in \cite{rudin}. (For general
spaces, the requirement of inner regularity is restricted to open sets and
sets of finite measure; see Theorem 2.14 in \cite{rudin}.)

Given $\lambda\geq1$ and $a\in X$, we say that a set $S_{\lambda}(a)\subseteq
X$ is a \textbf{Morse set} associated with $a$ and $\lambda$ if there is an
$r>0$ such that $B(a,r)\subseteq S_{\lambda}(a)\subseteq B(a,\lambda r)$ and
$S_{\lambda}(a)$ is starlike with respect to $B(a,r)$. We also say that
$S_{\lambda}(a)$ is a $\lambda$\textbf{-Morse set}. Recall that a gauge
function is a mapping $\delta\!:\!X\rightarrow(0,R)$ for some $R>0$. We will
say that the Morse set $S_{\lambda}(a)$ is $\delta$\textbf{-fine} with respect
to a gauge function $\delta$ if $\lambda r\leq\delta(a)$; in this case, we
will also call $a$ the \textbf{tag} for $S_{\lambda}(a)$. Note that putting
$\lambda=1$ forces a Morse set to be a closed ball. Also note that the closure
$\operatorname*{cl}(S_{\lambda}(a))$ of a $\lambda$-Morse set $S_{\lambda}(a)$
is again a $\lambda$-Morse set since when $y\in B(a,r)$, $x\in
\operatorname*{cl}(S_{\lambda}(a))$ and $\left\{  x_{n}\right\}  $ is a
sequence converging to $x$, we have for any $\alpha\in\lbrack0,1]$, $\alpha
y+\left(  1-\alpha\right)  x_{n}\rightarrow\alpha y+\left(  1-\alpha\right)
x$.

A collection $\mathcal{S}\subseteq\mathcal{P}(X)$ consisting of at least one
Morse set associated with each point $a$ in a set $\Omega\subseteq X$ is
called a \textbf{Morse cover} of $\Omega$ provided the same $\lambda\geq1$ is
used for each set in the cover and there is a finite upper bound to the
diameters of the sets in the cover. We will also call such a cover a $\lambda
$\textbf{-Morse cover.} A $\lambda$-Morse cover $\mathcal{S}$ of $\Omega$ is
called \textbf{fine} if for each $a\in\Omega$ and arbitrarily small values of
$r>0$ there are associated sets $S_{\lambda}(a)\in\mathcal{S}$ with
$B(a,r)\subseteq S_{\lambda}(a)\subseteq B(a,\lambda r)$ such that
$S_{\lambda}(a)$ is starlike with respect to $B(a,r)$. Given a Radon measure
$\mu$, a $\lambda$-Morse cover of a measurable set $\Omega\subseteq X$ is
called a $\mu$\textbf{-a.e.\hspace{0.03in}cover} of $\Omega$ if \ \textbf{i)}
it is fine, \ \textbf{ii)}\ each set in the cover is $\mu$-measurable, and
\ \textbf{iii)}\ for any $\varepsilon>0$, and any strictly positive gauge
function $\delta\!$ there is a finite or infinite sequence of disjoint,
$\delta$-fine sets $S_{n}\in\mathcal{S}$ such that $\mu(\Omega\setminus
\cup_{n}S_{n})=0$ and $\mu(\cup_{n}S_{n}\setminus\Omega)<\varepsilon$. This
concept is similar to that of Vitali covers, see \cite{federer}.

We first extend Corollary~\ref{finiteborelmeasure} to show that a fine Morse
cover consisting of closed sets is a $\mu$-a.e.\hspace{0.02in}cover for any
given measurable subset $\Omega$ of $X$. The same is true when the Morse sets
are not necessarily closed provided that for each set $E$ in the cover, it
does not increase the measure of $E$ to adjoin its closure points. For closed
balls and sets of finite measure, the proof is standard (see \cite{evans} or
\cite{ziemer}). We reproduce and extend it here.

\begin{lemma}
\label{lemma1} Let $\mu$ be a Radon measure on $X$. Let $\Omega\subseteq X$ be
measurable, and suppose that $\mathcal{S}$ is a fine Morse cover of $\Omega$
consisting of $\mu$-measurable sets. Then $\mathcal{S}$ is a $\mu$
-a.e.\hspace{0.02in}cover of $\Omega$ if $\mathcal{S}$ consists of closed sets
or if for each set $E\in\mathcal{S}$, $\mu(\Omega\cap(\operatorname*{cl}
(E)\setminus E))=0$.
\end{lemma}

\textbf{Proof:} Fix $\varepsilon>0$, and a gauge function $\delta>0$. We
suppose first that $\mathcal{S}$ consists of closed, $\delta$-fine sets. If
$\mu(\Omega)<\infty$, we may fix an open set $O\supseteq\Omega$ such that
$\mu(O\setminus\Omega)<\varepsilon$, and we may assume that each set
$E\in\mathcal{S}$ is a subset of $O$. Let $\kappa$ be the upper bound for the
Morse Covering Theorem; recall that it depends only on $X$ and the parameter
$\lambda$ for the cover. By Corollary~\ref{finiteborelmeasure}, there is a
finite subcollection $\mathcal{F}_{1}\subset\mathcal{S}$ consisting of
pairwise disjoint closed sets such that $\mu(\cup\mathcal{F}_{1})\geq
\mu(\Omega)/(2\kappa)$, whence $\mu(\Omega\setminus\cup\mathcal{F}_{1}
)\leq(1-1/(2\kappa))\mu(\Omega)$. Let $\Omega^{\prime}=\Omega\setminus
\cup\mathcal{F}_{1}$ and $\mathcal{S}_{1}=\{E\in\mathcal{S}:E\cap\left(
\cup\mathcal{F}_{1}\right)  =\varnothing\}$. Then $\mathcal{S}_{1}$ is a fine
Morse cover of $\Omega^{\prime}$. Again, there is a finite disjoint subfamily
$\mathcal{F}_{2}\subseteq\mathcal{S}_{1}$ such that $\mu(\Omega^{\prime
}\setminus\cup\mathcal{F}_{2})\leq(1-1/(2\kappa))\mu(\Omega^{\prime})$,
whence, $\mu(\Omega\setminus\cup(\mathcal{F}_{1}\cup\mathcal{F}_{2}
))\leq(1-1/(2\kappa))^{2}\mu(\Omega)$. Continuing in this manner, we have
$\mu(\Omega\setminus\mathcal{F})=0$ where $\mathcal{F}=\cup_{i}\mathcal{F}
_{i}$. Important for the next step, however, is the fact that for any
$\gamma>0$, there is a finite, pairwise disjoint family $\mathcal{F}^{\prime
}\subseteq\mathcal{S}$ such that $\mu(\Omega\setminus\mathcal{F}^{\prime
})<\gamma$.

Now suppose that $\mu(\Omega)=+\infty$. Then since $\mu$ is a Radon measure,
$\ \Omega=\cup_{i=1}^{\infty}\Omega_{i}$ where each $\Omega_{i}$ is a set of
finite measure and $\Omega_{i}\cap\Omega_{j}=\varnothing$ for $i\neq j$. For
each $i$, fix an open set $O_{i}\supseteq\Omega_{i}$ with $\mu\left(
O_{i}\setminus\Omega_{i}\right)  <\varepsilon/2^{i}$. We apply the above
result to obtain a finite (or empty) family $\mathcal{F}^{1}$ covering all but
a set of measure $1$ of $\Omega_{1}$ with all sets contained in $O_{1}$. At
the $n^{\text{th}}$ stage, $n>1$, we obtain a finite (or empty) family
$\mathcal{F}^{n}$ covering all but a set of measure $1/n$ of $\left(
\cup_{i=1}^{n}\Omega_{i}\right)  \setminus\cup_{i=1}^{n-1}\left(
\cup\mathcal{F}^{i}\right)  $ with all sets contained in $\left(  \cup
_{i=1}^{n}O_{i}\right)  \backslash\cup_{i=1}^{n-1}\left(  \cup\mathcal{F}
^{i}\right)  $. Clearly, $\cup_{i=1}^{\infty}\mathcal{F}^{i}$ is the desired
collection of disjoint sets in $\mathcal{S}$.

In the case that for each set $E\in\mathcal{S}$ , $\mu(\Omega\cap
(\operatorname*{cl}(E)\setminus E))=0$, we apply the above result to the Morse
cover formed by the closures of the sets in $\mathcal{S}$. We then replace
each set $\operatorname*{cl}(S_{n})$ in the resulting disjoint sequence with
the original set $S_{n}$. $\Box$

\bigskip

When dealing with Morse sets that are not closed, the conditions in
Lemma~\ref{lemma1} are easily fulfilled when the Morse cover $\mathcal{S}$ is
\textbf{scaled}. This means that for each $S_{\lambda}(a)\in\mathcal{S}$ and
each $p\in(0,1]$, the set $S_{\lambda}^{(p)}(a)$ is also in $\mathcal{S}$
where $S_{\lambda}^{(p)}(a)=\{a+px:a+x\in S_{\lambda}(a)\}$.

\begin{proposition}
\label{lemma2} Let $\mu$ be a Radon measure on $X$. Let $\Omega$ be a
measurable subset of $X$ and suppose $\mathcal{S}$ is a scaled Morse cover of
$\Omega$ consisting of $\mu$-measurable sets. Then $\mathcal{S}$ is a $\mu
$-a.e.\hspace{0.02in}cover of $\Omega$.
\end{proposition}

\textbf{Proof:} Since $\mathcal{S}$\ is a scaled Morse cover of $\Omega$, it
is certainly a fine cover of $\Omega$. Let $\lambda$ be the parameter for the
Morse cover $\mathcal{S}$. Let $a\in\Omega$ and fix $S_{\lambda}
(a)\in\mathcal{S}$; we write $S$ for $S_{\lambda}(a)$. We will show that for
$0<p<q\leq1$, $\partial S^{(p)}\cap\partial S^{(q)}=\varnothing$. The result
will then follow since for all but a countable number of values $p$,
$\mu(\partial S^{(p)})=0$. Since $S^{(p)}=\left(  S^{(q)}\right)  ^{(p/q)}$,
we may simplify notation by assuming that $S^{(q)}=S$; we may further simplify
by translating so that $a=0$. The result now follows from Proposition
\ref{cone} since for each $x\in\operatorname*{cl}(S^{(p)})$, $(1/p)x\in
\operatorname*{cl}(S)$, so $x\in\operatorname*{int}(S)$. $\Box$

\begin{example}
\begin{rm}
Take all closed balls or all open balls in $X$ of radius at most $1$. For each
center $x$ and radius $r$, let $a(x,r)$ in the interior of the ball be the tag
of that ball, and set $\omega(x,r):=\Vert x-a(x,r)\Vert/r$. Assume that
$\omega_{0}=\sup_{x,r}\omega(x,r)<1$. Given a Radon measure $\mu$, we have a
$\mu$-a.e.\hspace{0.02in}cover of any $\mu$-measurable set in $X$, and
$(1+\omega_{0})/(1-\omega_{0})$ is the smallest permissible value of $\lambda
$. As a special case, we may take each tag $a(x,r)=x$.
\end{rm}
\end{example}

\begin{example}
\begin{rm}
Let $\{e_{1},\ldots,e_{d}\}$ be a basis for $X$. Let $a=\sum_{i=1}^{d}
a_{i}e_{i}\in X$. Let $b$, $c\in\mathbb{R}_{+}^{d}=\{(x_{1},\ldots,x_{d}
)\in\mathbb{R}^{d}:x_{i}>0,1\leq i\leq d\}$ with $c_{1}^{2}+\cdots+c_{d}
^{2}<1$. Define a tagged interval by setting $I(a,b,c):=\{\sum_{i=1}^{d}
(a_{i}+t_{i})e_{i}:0<t_{i}\leq b_{i},1\leq i\leq d\}$ with tag at $\sum
_{i=1}^{d}(a_{i}+b_{i}c_{i})e_{i}$. Fix $c$ as above and take $k\geq1$. Given
a Radon measure $\mu$, the collection $\mathcal{S}=\{I(a,b,c):a\in
X,b\in\mathbb{R}_{+}^{d}$\ such that $\;\max_{1\leq i\leq d}b_{i}/\min_{1\leq
i\leq d}b_{i}\leq k\}$ is a scaled, $\mu$-a.e.\hspace{0.03in}\hspace
{0.02in}Morse cover of $X$.
\end{rm}
\end{example}

Let $K$ be a compact subset of $X$ and let $\mu$ be a Radon measure such that
each open ball with center at a point of $K$ has positive $\mu$-measure. We
will want to use the fact that given a $\lambda\geq1$, any $\mu$-a.e.,
$\lambda$-Morse cover $\mathcal{S}$ of $K$ forms a differentiation basis on
$K$ with respect to $\mu$. For our purposes here this means that if $\nu$ is a
radon measure absolutely continuous with respect to $\mu$, i.e., $\nu<<\mu$,
and $\mathcal{S}(a)$ is the collection of sets in $\mathcal{S}$ associated
with $a\in K$, then
\[
\lim_{\substack{S\in\mathcal{S}(a)\\\Delta(S)\rightarrow0}}\frac{\nu(S)}
{\mu(S)}=\frac{d\nu}{d\mu}(a)\;\;\text{for \ \ }\mu\text{-a.e.\ \ }a\in K,
\]
where $\frac{d\nu}{d\mu}$ denotes the Radon-Nikod\'{y}m derivative of $\nu$
with respect to $\mu$.

By the principal result in \cite{bliedtnerloeb1}, the above equality follows
from the fact that if $E$ is a measurable subset of $K$ and $\nu$ is a finite
Radon measure with $\nu<<\mu$ and $\nu(E)=0$, then for $\mu$-a.e.$\hspace
{0.04in}a\in E$, $\lim\sup_{S\in S(a),\,\Delta(S)\rightarrow0}\nu
(S)/\mu(S)\leq1$. As in \cite{bliedtnerloeb1}, we can see that this is in fact
the case by letting $A$ be the subset of $E$ where the reverse inequality
holds, and letting $\kappa$ be the upper bound given by the Morse Covering
Theorem. We fix $\mathbb{\varepsilon}>0$ and a nonempty compact set $C\subset
X\setminus E$ with $\nu(X\setminus C)<\mathbb{\varepsilon}/\left(
2\kappa\right)  $. By assumption, for each $a\in A$, there is a set
$S(a)\in\mathcal{S}(a)$ with $S(a)\cap C=\varnothing$ and $\mu(S(a))\leq
\nu(S(a))$. For the finite, disjoint subcollection $\left\langle
S_{n}\right\rangle $ of these sets given by Corollary 
\ref{finiteborelmeasure}
, we have
\[
\mu^{\ast}(A)\leq2\kappa\cdot\Sigma_{n}\mu(S_{n})\leq2\kappa\cdot\Sigma_{n}
\nu(S_{n})\leq2\kappa\cdot\nu(X\setminus C)<\mathbb{\varepsilon}\text{.}
\]

In the next section, we will want to exploit the fact that measurable
functions are approximately continuous almost everywhere with respect to a
given Radon measure $\mu$. That is, let $\Omega$ be a $\mu$-measurable subset
of $X$, and let $f\!:\!\Omega\rightarrow\mathbb{R}$ be $\mu$-measurable; set
$f\equiv0$ on $X\setminus\Omega$. Suppose $\mathcal{S}$ is a fine $\lambda
$-Morse cover of $\Omega$, so that the sets in $\mathcal{S}$ form a
differentiation basis with respect to $\mu$ at points $x\in\Omega$ for which
all balls $B(x,r)$ have positive $\mu$-measure. Then $x\in\Omega$ is called a
\textbf{point of approximate continuity} of $f$ if for all positive
$\varepsilon$ and $\eta$ there is an $R>0$ such that if $S(x)$ is a set in
$\mathcal{S}$ with tag $x$ and $S(x)\subseteq B(x,R)$, then for $E(x,\eta
):=\{y\in S(x):|f(x)-f(y)|>\eta\}$ we have $\mu(E(x,\eta))\leq\varepsilon
\,\mu(S(x))$. It is known that $\mu$-almost all points of $\Omega$ are points
of approximate continuity of $f$ (see \cite{federer}, 2.9.13). A related
notion, defined and used below in the proof of Theorem \ref{integration}, is
the notion of a \textbf{Lebesgue point} for $f$; these also fill the space
except for a set of measure $0$.

\begin{remark}
\begin{rm}
Clearly, a nonnegative, measurable, real-valued function $f$ is approximately
continuous $\mu$-a.e.\hspace{0.02in}if for each $n\in\mathbb{N}$, $\min\left(
f,n+1\right)  $ is approximately continuous $\mu$-a.e.\hspace{0.02in}on the
set where $f\leq n$. That this is the case follows from the discussion of
Lebesgue points in Section 3 of \cite{bliedtnerloeb2}, since the constant for
a Lebesgue point $x$ equals $f(x)$ for $\mu$-almost all $x$ (cf. Equation
\eqref{2} below).
\end{rm}
\end{remark}

\section{Integration}

\label{Integration}

Again, we let $(X,\Vert\cdot\Vert)$ be a normed vector space of dimension
$d<\infty$ over the real numbers $\mathbb{R}$. Using our covering results we
can formulate the Lebesgue integral as a type of Riemann sum defined by $\mu
$-a.e.\hspace{0.02in}Morse covers. We do this first for nonnegative functions
and later apply the result to measurable functions taking both positive and
negative values.

\begin{theorem}
\label{integration} Let $\mu$ be a Radon measure on $X$. Let $\Omega$ be a
measurable subset of $X$, and let $f$ be a nonnegative, real-valued,
measurable function on $\Omega$. Then $\int_{\Omega}f\,d\mu$ is finite and
equals $F$ if the following condition holds for some $\lambda\geq1$ and some
$\mu$-a.e., $\lambda$-Morse cover $\mathcal{S}$ of $\Omega\!:$ For all
$\varepsilon>0$ there is a gauge function $\delta\!:\!\Omega\rightarrow(0,1]$
such that for any finite or countably infinite disjoint sequence $\left\langle
S_{n}(x_{n})\right\rangle $ of $\delta$-fine sets from $\mathcal{S}$ covering
all but a set of measure $0$ of $\Omega$ we have
\begin{equation}
\left|  \sum\nolimits_{n}\,f(x_{n})\,\mu(S_{n})-F\right|  <\varepsilon.
\label{1}
\end{equation}
Conversely, if $\int_{\Omega}f\,d\mu$ is finite and equals $F$, then the
condition holds for any $\lambda\geq1$ and any $\mu$-a.e., $\lambda$-Morse
cover $\mathcal{S}$ of $\Omega$.
\end{theorem}

\textbf{Proof:} We note first that\ for a given set $A\subseteq\Omega$ with
$\mu(A)=0$, we may set our gauge to force an arbitrarily small sum for points
$x_{i}\in A$, and also force, in the case that $f$ is assumed to be
integrable, an arbitrarily small integral of $f$ over the union of the
corresponding sets $S_{i}$. To show this, we fix $\varepsilon>0$, and for each
$n\in\mathbb{N}$ we set $A_{n}=\{x\in A:n-1\leq f(x)<n\}$. The sets $A_{n}$
are disjoint and $\mu$-null with union $A$. In the case that $f$ is assumed to
be integrable, we may choose an open set $G\supseteq A$ so that $\int
_{G}f<\varepsilon$; otherwise, set $G=X$. For each $n\in\mathbb{N}$, fix an
open set $G_{n}$ with $G\supseteq$ $G_{n}\supseteq A_{n}$ and $\mu
(G_{n})<\varepsilon/\left(  n\cdot2^{n}\right)  $. (This is possible since
$\mu$ is outer regular.) For each $x\in A_{n}$, we choose $\delta
(x)<\sup\{s:B(x,s)\subseteq G_{n}\}$. Then a sum over $\delta$-fine, disjoint
sets $S_{i}$ with all tags in $A$ satisfies the inequality
\[
\sum\limits_{i}f(x_{i})\,\mu(S_{i})<\sum\limits_{n=1}^{\infty}\left(
n\sum\limits_{x_{i}\in A_{n}}\!\!\mu(S_{i})\right)  \leq\sum\limits_{n=1}
^{\infty}\varepsilon\,2^{-n}=\varepsilon,
\]
and if $f$ is assumed to be integrable, its integral over $\cup_{i}S_{i}$ is
at most $\varepsilon$.

Now suppose that $\int_{\Omega}f\,d\mu$ exists and equals $F$. Fix
$\lambda\geq1$, a $\mu$-a.e., $\lambda$-Morse cover $\mathcal{S}$ of $\Omega$,
and an $\varepsilon>0$. Set $f\equiv0$ on $X\setminus\Omega$. Suppose
$x\in\Omega$ \ is a Lebesgue point for $f$ with respect to the Morse cover
$\mathcal{S}$. This means that there is a constant, which (after redefining
$f$ on a $\mu$-null set) we may assume is $f(x)$, such that the following
condition holds: For any $\varepsilon_{1}>0$ there is an $R>0$ so that
if\ $S(x)$ is a set in $\mathcal{S}$ with tag $x$ and $S(x)\subseteq B(x,R)$,
then
\begin{equation}
\int\limits_{y\in S(x)}\!\!|f(x)-f(y)|\,d\mu(y)\leq\varepsilon_{1}\,\mu(S(x)).
\label{2}
\end{equation}
For such an $x$, let $k(x)$ be the first integer strictly larger than
$\left\|  x\right\|  $. Set $\delta(x)=R$ where $R$ is chosen to be at most
$1$ and satisfy Equation \eqref{2} with $\varepsilon_{1}=\varepsilon
\,2^{-k(x)-1}\,/\left[  1+\mu(B(0,k(x)+1))\right]  $. Since $\mathcal{S}$
forms a differentiation basis, it follows that the non-Lebesgue points form a
$\mu$-null set. (See, for example, Section 3 of \cite{bliedtnerloeb2}.) We
may, as just noted, choose positive values $\delta(x)\leq1$ for such points
$x$ so that their contribution to the sum in Equation \eqref{1} can be at most
$\varepsilon/4$ and the integral of $f$ over the union of the corresponding
sets $S(x)$ will be at most $\varepsilon/4$.

With this choice for the gauge $\delta$, we now let $\left\langle S_{n}
(x_{n})\right\rangle $ be any finite or countably infinite disjoint sequence
of $\delta$-fine sets from $\mathcal{S}$ covering all but a set of measure $0$
of $\Omega$. Let $L$ denote the set of Lebesgue points of $\Omega$. Then

\begin{align}
\left|  \int_{\Omega}f\,d\mu-\sum\limits_{n}f(x_{n})\,\mu(S_{n})\right|   &
=\left|  \,\,\int\limits_{\cup_{n} S_{n}}\!\!f\,d\mu-\sum\limits_{n}
f(x_{n})\,\mu(S_{n})\right| \label{3}\\
&  \leq\sum\limits_{x_{n}\in L}\int_{S_{n}}\!\!|f(x_{n})-f(y)|\,d\mu
(y)+\frac{\varepsilon}{2}\nonumber\\
&  \leq\sum\limits_{\ell=1}^{\infty}\frac{\varepsilon\,2^{-\ell-1}}
{1+\mu(B(0,\ell+1))}\sum\limits_{\ell-1\leq\left\|  x_{n}\right\|  <\ell}
\mu(S_{n})+\frac{\varepsilon}{2}\nonumber\\
&  \leq\varepsilon.\nonumber
\end{align}

Now fix a $\lambda\geq1$ and a $\mu$-a.e., $\lambda$-Morse cover $\mathcal{S}$
of $\Omega$ so that for any $\varepsilon>0$ there is an appropriate gauge
$\delta\leq1$ for $f$ and $F$; that is, for any finite or countably infinite
disjoint sequence $\left\langle S_{n}(x_{n})\right\rangle $ of $\delta$-fine
sets from $\mathcal{S}$ covering all but a set of measure $0$ of $\Omega$,
Equation \eqref{1} holds for $\varepsilon$. For each $x\in\Omega$, let $k(x)$
be the first integer strictly larger than $\left\|  x\right\|  $, and set
\[
\eta(x):=\frac{2^{-k(x)}}{[1+\mu(B(0,k(x)+1))][1+f(x)]}.
\]
For each $m\in\mathbb{N}$, fix $\delta_{m}\leq1$ to work for $f$ and $F$ with
$\varepsilon=1/m$ in Equation \eqref{1}. Let $\left\langle S_{n}^{m}(x_{n}
^{m})\right\rangle $ be a finite or countably infinite disjoint sequence of
$\delta_{m}$-fine sets from $\mathcal{S}$ covering all but a set of measure
$0$ of $\Omega$. We may assume that each tag $x_{n}^{m}$ is a point of
approximate continuity of $f$ \ and $\mu(E_{n}^{m})\leq\eta_{m}(x_{n}
^{m})\,\mu(S_{n}^{m})$ where $\eta_{m}(x_{n}^{m})=\eta(x_{n}^{m})/m$ and
\[
E_{n}^{m}:=\{x\in S_{n}^{m}:|f(x_{n}^{m})-f(x)|>\eta_{m}(x_{n}^{m})\}\text{.}
\]
Define a measurable function $f_{m}$ on $\Omega$ as follows: If for some
$n\in\mathbb{N}$, $x\in S_{n}^{m}\setminus E_{n}^{m}$, set $f_{m}
(x)=\max(f(x_{n}^{m})-\eta_{m}(x_{n}^{m}),0)$; otherwise, set $f_{m}(x)=0$.
Now the functions $f_{m}$ converge to $f$ in measure since,
\begin{align*}
&  \mu\left(  \left\{  x\in\Omega:|f(x)-f_{m}(x)|>\frac{1}{m}\right\}  \right)
\\
&  \leq\sum\limits_{n}\mu(E_{n}^{m})\leq\sum\limits_{n}\eta_{m}(x_{n}
^{m})\,\mu(S_{n}^{m})\\
&  \leq\frac{1}{m}\sum\limits_{\ell=1}^{\infty}\frac{2^{-\ell}}{1+\mu
(B(0,\ell+1))}\sum\limits_{\ell-1\leq\left\|  x_{n}^{m}\right\|  <\ell
}\!\!\!\!\!\mu(S_{n}^{m})\\
&  \leq\frac{1}{m}.
\end{align*}
Since any subsequence of the sequence $\left\langle f_{m}\right\rangle $ has
in turn a subsequence converging $\mu$-a.e.\hspace{0.02in}to $f$, it follows
from Fatou's lemma that
\begin{align*}
\int_{\Omega}f\,d\mu &  \leq\lim\inf_{m}\int_{\Omega}f_{m}\,d\mu\\
&  \leq\lim\inf_{m}\sum\limits_{n}f(x_{n}^{m})\,\mu(S_{n}^{m})\\
&  \leq\lim\inf_{m}\,(F+1/m)\\
&  =F<+\infty.
\end{align*}
On the other hand, each $f_{m}\leq f$ , so for each $m$,
\begin{align*}
\int_{\Omega}f\,d\mu &  \geq\int_{\Omega}f_{m}\,d\mu\\
&  \geq\sum\limits_{n}\left[  f(x_{n}^{m})-\eta_{m}(x_{n}^{m})\right]
\,\mu(S_{n}^{m}\setminus E_{n}^{m})\\
&  =\sum\limits_{n}f(x_{n}^{m})\,\mu(S_{n}^{m})-\sum\limits_{n}f(x_{n}
^{m})\,\mu(E_{n}^{m})-\sum\limits_{n}\eta_{m}(x_{n}^{m})\,\mu(S_{n}
^{m}\setminus E_{n}^{m})\\
&  \geq F-1/m-\sum\limits_{n}\eta_{m}(x_{n}^{m})\,f(x_{n}^{m})\,\mu(S_{n}
^{m})-\sum\limits_{n}\eta_{m}(x_{n}^{m})\,\mu(S_{n}^{m})\\
&  \geq F-1/m-\frac{1}{m}\sum\limits_{\ell}\frac{2^{-\ell}}{1+\mu
(B(0,\ell+1))}\sum\limits_{\ell-1\leq\left\|  x_{n}^{m}\right\|  <\ell}
2\,\mu(S_{n}^{m})\\
&  \geq F-3/m,
\end{align*}
whence $\int_{\Omega}f\,d\mu=F$. $\Box$

\begin{remark}
\begin{rm}
With no loss of generality, we can restrict the points $x_{n}$ in Equation
\eqref{1} to be points of approximate continuity or to be points outside of
any given $\mu$-null set. Also, while we could work with the cover formed by
all $\mu$-measurable $\lambda$-Morse sets, \ the gauge $\delta$ can in general
be chosen larger when given a smaller $\mu$-a.e.\hspace{0.03in}Morse cover.
\end{rm}
\end{remark}

Let $f\!$ be a real-valued function on $\Omega$ taking both positive and
negative values. As usual, we set $f^{+}:=\max(f,0)$ and $f^{-}:=\max(-f,0)$;
given $\mu$, we say that $f$ is integrable if both $f^{+}$ and $f^{-}$ have
finite integrals with respect to $\mu$. Suppose now that $\mathcal{S}$ is the
set of all closed balls in $X$ with tags at the center; i.e., $\lambda=1$.
Even for this case, we cannot force the integrability of $f$ with the
inequality
\begin{equation}
\left|  \sum\limits_{n=1}^{\infty}f(x_{n})\,\mu(B_{n})-F\right|  <\varepsilon.
\label{4}
\end{equation}
The inequality does imply that $\sum|f(x_{n})|\mu(B_{n})$ will be bounded for
any appropriate sequence $\left\langle B_{n}\right\rangle $, but the sums need
not be uniformly bounded. The condition given by \eqref{4} will allow
principal value integrals. For example, in $\mathbb{R}^{d}$ let $e_{1}$ be the
unit vector in the positive direction along the first axis. For each
$n\in\mathbb{N}$, let $A_{n}$ be the open ball $U\left(  \left(
(-1)^{n}/n\right)  \cdot e_{1},1/(2n^{2})\right)  $. The balls $A_{n}$ are
disjoint. Let $\Omega$ be the union of the balls $A_{n}$ together with the
origin, and let $\mu$ be the sum of the Dirac measure supported at the origin
and Lebesgue measure restricted to $\Omega$. Set $f(x)=\left(  (-1)^{n}
/n\right)  \cdot\mu\left(  A_{n}\right)  $ if $x\in A_{n}$, and let $f(0)=0$.
Take the gauge function $\delta\!:\!\mathbb{R}^{d}\rightarrow(0,1)$ so that if
$x\in A_{n}$ then $B(x,\delta(x))\subset A_{n}$. Let $F=\sum_{n=1}^{\infty
}(-1)^{n}/n$, i.e., $F=-\ln2$. Given $\varepsilon>0$, if we take $\delta$
small enough at $0$, then for any sequence of disjoint, $\delta$-fine balls
$B_{n}$ satisfying $\mu(\Omega/\cup_{n}B_{n})=0$, we have $|\sum f(x_{n}
)\,\mu(B_{n})-F|<\varepsilon$. Any such sequence must contain a ball having
the origin as its center. As we choose different sequences so that the radius
of this ball shrinks to $0$ we have $\sum|f(x_{n})|\mu(B_{n})\rightarrow
\infty$.

It is the case, as we now show, that a real-valued, measurable $f$ is
integrable if the sums $\sum\nolimits_{n}\,\left|  f(x_{n})\right|
\,\mu(S_{n})$ are uniformly bounded.

\begin{corollary}
\label{UniformBound}Given $\mu$ and $\Omega$ as in the theorem, let $f\!$ be
an arbitrary, measurable, real-valued function on $\Omega$. Then $f$ is
integrable if the following condition holds for some $\lambda\geq1$ and some
$\mu$-a.e., $\lambda$-Morse cover $\mathcal{S}$ of $\Omega\!:$ There is a
number $M\geq0$ and a gauge function $\delta\!:\!\Omega\rightarrow(0,1]$ such
that for any finite or countably infinite disjoint sequence $\left\langle
S_{n}(x_{n})\right\rangle $ of $\delta$-fine sets from $\mathcal{S}$ covering
all but a set of measure $0$ of $\Omega$ we have
\begin{equation}
\sum\nolimits_{n}\,\left|  f(x_{n})\right|  \,\mu(S_{n})\leq M\text{.}
\label{5}
\end{equation}
Conversely, if $f$ is integrable, then the condition holds for all
$\lambda\geq1$ and all $\mu$-a.e., $\lambda$-Morse covers $\mathcal{S}$ of
$\Omega$. In this case, \ for each such $\lambda$-Morse cover $\mathcal{S}$
and each $\varepsilon>0$, there is a gauge function $\delta\!:\!\Omega
\rightarrow(0,1]$ so that for any finite or countably infinite disjoint
sequence $\left\langle S_{n}(x_{n})\right\rangle $ of $\delta$-fine sets from
$\mathcal{S}$ covering all but a set of measure $0$ of $\Omega$ we have
\[
\left|  \sum\nolimits_{n}\,f(x_{n})\,\mu(S_{n})-\int_{\Omega}f\,d\mu\right|
<\varepsilon.
\]
\end{corollary}

\textbf{Proof:} Fix a $\lambda\geq1$ and a $\mu$-a.e., $\lambda$-Morse cover
$\mathcal{S}$ of $\Omega$, and suppose there is an $M\geq0$ and a gauge
$\delta\leq1$ satisfying our condition including Equation \eqref{5}. For each
$x\in\Omega$, let $k(x)$ be the first integer strictly larger than $\left\|
x\right\|  $, and set $\eta(x):=2^{-k(x)}/\left(  1+\mu(B(0,k(x)+1))\right)
$. For each $m\in\mathbb{N}$, let $\left\langle S_{n}^{m}(x_{n}^{m}
)\right\rangle $ be a finite or countably infinite disjoint sequence of
$\delta$-fine sets from $\mathcal{S}$ covering all but a set of measure $0$ of
$\Omega$. We may assume that each tag $x_{n}^{m}$ is a point of approximate
continuity of $\left|  f\right|  $ \ and $\mu(E_{n}^{m})\leq\eta_{m}(x_{n}
^{m})\,\mu(S_{n}^{m})$ where $\eta_{m}(x_{n}^{m})=\eta(x_{n}^{m})/m$ and
\[
E_{n}^{m}:=\{x\in S_{n}^{m}:\left|  \left|  f(x_{n}^{m})\right|  -\left|
f(x)\right|  \right|  >\eta_{m}(x_{n}^{m})\}\text{.}
\]
Define a measurable function $f_{m}$ on $\Omega$ as follows: If for some
$n\in\mathbb{N}$, $x\in S_{n}^{m}\setminus E_{n}^{m}$, set $f_{m}
(x)=\max(\left|  f(x_{n}^{m})\right|  -\eta_{m}(x_{n}^{m}),0)$; otherwise, set
$f_{m}(x)=0$. As in the theorem, we have $f_{m}\rightarrow|f|$ in measure and
\[
\int_{\Omega}|f|\,d\mu\leq\lim\inf\int_{\Omega}f_{m}\,d\mu\leq\lim\inf\sum
_{n}|f(x_{n}^{m})|\,\mu(S_{n}^{m})\leq M,
\]
whence $f$ is integrable.

Now assume that $f$ is integrable, and set $F_{1}=\int_{\Omega}|f|\,d\mu$.
Applying the theorem, it follows that for any $\lambda\geq1$ and any $\mu
$-a.e., $\lambda$-Morse cover $\mathcal{S}$ of $\Omega$, the function $\left|
f\right|  $ satisfies our condition including Equation \eqref{5} with
$M=F_{1}+1$. The rest follows for any given $\varepsilon>0$ by applying the
theorem separately to $f^{+}$and $f^{-}$ with respect to $\varepsilon/2$ and
taking the smaller of the two gauges at each point. $\Box$

\section{An Extension of the Riemann Integral}

\label{Riemann}

For the case that $f$ is real-valued and continuous almost everywhere, we can
easily calculate the gauge $\delta$, and in the process obtain an extension of
the Riemann integral that integrates some unbounded functions with respect to
Radon measures on unbounded domains. Here too, we say that $f$ is integrable
only when this is true for $f^{+}$ and $f^{-}$.

\begin{theorem}
Let $\mu$ be a Radon measure on $X$. Let $\Omega$ be a measurable subset of
$X$, and let $f$ be a measurable, real-valued function on $\Omega$. Set
$f\equiv0$ on $X\setminus$ $\Omega$, and let $\Omega_{c}$ be the set of points
in $\Omega$ where $f$ is continuous. Let us suppose that $\mu(\Omega
\setminus\Omega_{c})=0$. For $x\in\Omega_{c}$, let $k(x)$ be the smallest
integer strictly greater than $\left\|  x\right\|  $, and for each $\gamma>0$
fix $\rho(x,\gamma)$ with $0<\rho(x,\gamma)\,\leq1$ so that for all $y$ with
$\left|  y-x\right|  <\rho(x,\gamma)$, we have $|f(y)-f(x)|<\gamma$. If
$\mu(\Omega)<\infty$, then for each $\varepsilon>0$ and each $x\in\Omega_{c}$
set $\delta_{\varepsilon}(x)=\rho(x,\varepsilon\cdot\,[1+\mu(\Omega)]^{-1})$;
otherwise for each $\varepsilon>0$ and each $x\in\Omega_{c}$ set
$\delta_{\varepsilon}(x)=\rho(x,\varepsilon\cdot2^{-k(x)}\cdot\lbrack
1+\mu(B(0,k(x)+1))]^{-1})$. Now, if $f$ is integrable, then for any
$\lambda\geq1$, any $\varepsilon$ with $0<\varepsilon\leq1$ and any finite or
countably infinite disjoint sequence $\left\langle S_{n}(x_{n})\right\rangle $
of $\delta_{\varepsilon}$-fine, $\lambda$-Morse sets covering all but a set of
measure $0$ of $\Omega$ and having tag points $x_{n}\in\Omega_{c}$ we have
\[
\left|  \sum\nolimits_{n}\,f^{\pm}(x_{n})\,\mu(S_{n}(x_{n}))-\int_{\Omega
}f^{\pm}\,d\mu\right|  <\varepsilon,
\]
whence
\[
\left|  \sum\nolimits_{n}\,f(x_{n})\,\mu(S_{n}(x_{n}))-\int_{\Omega}
f\,d\mu\right|  <2\varepsilon.
\]
On the other hand, $f$ is integrable if for some $\lambda\geq1$ and some
finite or countably infinite disjoint sequence $\left\langle S_{n}
(x_{n})\right\rangle $ of $\delta_{1}$-fine Morse sets, associated with
$\lambda$ and tag points $x_{n}\in\Omega_{c}$, and covering all but a set of
measure $0$ of $\Omega$, we have
\[
\sum\nolimits_{n}\,\left|  f(x_{n})\right|  \,\mu(S_{n})<+\infty.
\]
\end{theorem}

\textbf{Proof:} Note that if $\rho(x,\gamma)$ works for $f$, then it works for
$f^{+}$ and $f^{-}$. Assume $f$ is integrable, and fix $\lambda\geq1$ and
$\varepsilon>0$. Let $\left\langle S_{n}(x_{n})\right\rangle $ be any finite
or countably infinite disjoint sequence of $\delta_{\varepsilon}$-fine,
$\lambda$-Morse sets with tag points $x_{n}$ in $\Omega_{c}$ and covering all
but a set of measure $0$ of $\Omega$. Then for the case that $\mu
(\Omega)=\infty$ we have
\begin{align*}
\left|  \int_{\Omega}f^{+}\,d\mu-\sum\limits_{n}f^{+}(x_{n})\,\mu
(S_{n})\right|   &  \leq\sum\limits_{n}\int_{S_{n}}\!\!|f^{+}(y)-f^{+}
(x_{n})|\,d\mu(y)\\
&  \leq\sum\limits_{\ell=1}^{\infty}\frac{\varepsilon\,2^{-\ell}}
{1+\mu(B(0,\ell+1))}\sum\limits_{\ell-1\leq\left\|  x_{n}\right\|  <\ell}
\mu(S_{n})\leq\varepsilon,
\end{align*}
with the obvious simplification for the case that $\mu(\Omega)<\infty$. A
similar calculation works for $f^{-}$.

Now fix $\lambda\geq1$, and assume there is a finite or countably infinite
disjoint sequence $\left\langle S_{n}(x_{n})\right\rangle $ of $\delta_{1}
$-fine Morse sets associated with $\lambda$ and tag points $x_{n}\in\Omega
_{c}$ covering all but a set of measure $0$ of $\Omega$ such that
$\sum\nolimits_{n}\,\left|  f(x_{n})\right|  \,\mu(S_{n})=M\in\mathbb{R}$.
Then for the case that $\mu(\Omega)=\infty$,
\begin{align*}
\int_{\Omega}\left|  f\right|  \,d\mu &  =\sum_{n}\int_{S_{n}}\left|
f\right|  \,d\mu\\
&  \leq\sum\limits_{\ell}\left(  \left|  f(x_{n})\right|  +\frac{\,2^{-\ell}
}{1+\mu(B(0,\ell+1))}\right)  \sum_{\ell-1\leq\left\|  x_{n}\right\|  <\ell
}\!\!\!\!\!\mu(S_{n})\\
&  \leq M+1\text{.}
\end{align*}
Again, we have the obvious simplification for the case that $\mu
(\Omega)<\infty$. $\Box$

\bigskip
\noindent
{\bf Note added in proof:}  It follows from Lusin's Theorem and the Lebesgue
Differentiation Theorem for characteristic functions that this theory
can be extended to Banach space valued functions.  This will be the
subject of a subsequent paper.


\begin{thebibliography}{99}
\bibitem{besicovitch}A.S. Besicovitch, \textit{A general form of the covering
principle and relative differentiation of additive functions} (I), (II),
{Proc. Cambridge Phil. Soc. \textbf{41} (1945) 103--110, \textbf{42} (1946)
1--10. }

\bibitem {bliedtnerloeb1}J. Bliedtner and P. Loeb, \textit{A reduction
technique for limit theorems in analysis and probability theory}, Ark. Mat.
\textbf{30} (1992) 25--43.

\bibitem {bliedtnerloeb2}J. Bliedtner and P. Loeb, The optimal differentiation
basis and liftings of $L^{{\infty}}$, to appear in \textit{Trans. Amer. Math.
Soc}.

\bibitem {evans}L.C. Evans and R.F. Gariepy, \textit{Measure theory and fine
properties of functions}, Boca Raton, CRC Press, 1992.

\bibitem {federer}H. Federer, \textit{Geometric measure theory}, Berlin,
Springer--Verlag, 1969.

\bibitem {folland}G.B. Folland, \textit{Real analysis}, New York, Wiley, 1999.

\bibitem {furediloeb}Z. F\"{u}redi and P.A. Loeb, \textit{On the best constant
for the Besicovitch covering theorem}, Proc. Amer. Math. Soc., \textbf{121}
(1994) 1063--1073.

\bibitem {gordon}R.A. Gordon, \textit{The integrals of Lebesgue, Denjoy,
Perron, and Henstock}, Providence, American Mathematical Society, 1994.

\bibitem {loeb}P.A. Loeb, \textit{Opening the covering theorems of Besicovitch
and Morse}, Mathematica Moravica, Special volume (1997) 3--11.

\bibitem {malee}Z.M. Ma and P.Y. Lee, \textit{Absolute integration using
Vitali covers}, Real Anal. Exchange \textbf{18} (1992/93) 409--419.

\bibitem {morse}A.P. Morse, \textit{Perfect blankets}, Trans. Amer. Math. Soc.
\textbf{61} (1947) 418--442.

\bibitem {rudin}W. Rudin, \textit{Real and complex analysis}, New York,
McGraw--Hill, 1987.

\bibitem {ziemer}W.P. Ziemer, \textit{Weakly differentiable functions}, New
York, Springer--Verlag, 1989.
\end{thebibliography}
\end{document}